\documentclass{fundam}
\usepackage{amssymb}
\usepackage{enumerate}
\usepackage{mymacros}

\begin{document}

\setcounter{page}{155}
\publyear{22}
\papernumber{2124}
\volume{186}
\issue{1-4}

            \finalVersionForARXIV

\title{On the Tutte and Matching Polynomials for Complete Graphs}

\author{Tomer Kotek\\
Berlin, Germany \\
 tomer.kotek@gmail.com
\and Johann A. Makowsky\thanks{Address for correspondence:  Department of Computer Science,
              Israel Institute of Technology, Haifa, Israel}
 \\
Department of Computer Science \\
Israel Institute of Technology, Haifa, Israel\\
 janos@cs.technion.ac.il
}

\maketitle

\runninghead{T. Kotek and J.A. Makowsky}{On the Tutte Polynomial}

\begin{abstract}
Let $T(G;X,Y)$ be the Tutte polynomial for graphs.
We study the sequence $t_{a,b}(n) = T(K_n;a,b)$ where $a,b$ are integers, and show that
for every $\mu \in \N$ the sequence $t_{a,b}(n)$
is ultimately periodic  modulo $\mu$ provided $a \neq 1 \mod{\mu}$ and $b \neq 1 \mod{\mu}$.
This result is related to a conjecture by A. Mani and R. Stones from 2016.
The theorem is a consequence of a more general theorem which holds for a wide class of
graph polynomials definable in Monadic Second Order Logic. This gives  also similar results for the
various substitution instances of the
bivariate matching polynomial and the trivariate edge elimination polynomial $\xi(G;X,Y,Z)$ introduced by
I. Averbouch, B. Godlin and the second author in 2008.
All our results depend on the Specker-Blatter Theorem from 1981, which studies modular recurrence
relations of combinatorial sequences which count the number of labeled graphs.
\end{abstract}

\section{Introduction}

Boris (Boaz) Abramovich Trakhtenbrot (1921-2016) was one of the pioneers
in recognizing the usefulness of Monadic Second Order Logic $\MSOL$
for treating situations in automata theory \cite{trakhtenbrot1962finite}.
In \cite{ar:FischerMakowsky08}, published in the Festschrift for Boaz' 85th birthday,
Eldar Fischer and the second author gave an application of Monadic Second Order Logic
to graph polynomials, by proving that for a  graph polynomial $P(G;\bar{x})$ definable in $\MSOL$
and a sequence of recursively defined graphs $G_i: i \in \N$
the sequence of polynomials $P(G_i;\bar{x})$ satisfies a linear recurrence relation over the
polynomial ring $\Z[\bar{x}]$.

Ernst Specker (born in 1920-2011), together with Christian Blatter (1935-2021), was the first to use
Monadic Second Order Logic to prove a meta-theorem in counting combinatorics, \cite{pr:BlatterSpecker84}.
In order to celebrate Boaz' centenary we shall give another application of Monadic Second Order Logic
to $\MSOL$-definable  graph polynomials $P(G;\bar{x})$, where the proof again uses the Specker-Blatter Theorem.
We fix $\bar{a} \in \N^r$ and look at the sequence  $p(n, \bar{a})=P(K_n;\bar{a})$ modulo $\mu \in \N$,
where $K_n$ is the complete
graph on $n$ vertices. We show that, under simple assumption on $\bar{a}$ and $\mu$,
this sequence is ultimately periodic modulo $\mu$.

Our results
are easy, but possibly unexpected,
applications of the Specker-Blatter Theorem.
They illustrate more the power of it,
by using general methods of Monadic Second Order Logic $\MSOL$, rather than
applying combinatorial arguments specially tailored to a particular case.

\subsection{Some graph polynomials}

Let $G =(V(G), E(G))$ be a finite graph.
We put $n(G) = |V(G)|$, $m(G) = |E(G)|$, $\kappa(G)$ is the number of connected components of $G$.
From a logical point of view $G$ can be represented in various ways.
The vocabulary $\tau_{graph}$ consists of one binary relation symbol for $E(G)$
and the vocabulary $\tau_{hgraph}$ consists of two unary predicates, one for vertices $V(G)$ and one for edges $E(G)$
and  one binary relation symbol for $R(G)$ for the incidence relation between edges and vertices.
$\tau_{hgraph}$ is also suitable for representing hypergraphs.
{\em Induced subgraphs} are substructures in the case of $\tau_{graph}$.
and {\em subgraphs} are substructures in the case of $\tau_{hgraph}$.
In other words, let $A \subseteq V(G)$ be a set of vertices.
The {\em induced subgraph} of $G$
generated by $A$ is the graph
$G[A] = (A, E(G) \cap A^2)$.
On the other hand the graph $(A,E)$ is a subgraph of $G$ for any
$F \subseteq E(G)$.

$\MSOL$ on graphs allows quantification over subsets of vertices only.
$\MSOL$ on hypergraphs allows quantification over subsets of vertices and subsets of edges.
$\MSOL$ on hypergraphs has the same expressive power as Second Order Logic $\SOL$ on graphs where second order quantification
is {\em restricted to unary predicates and binary relations which are subsets of the edge relation}.
We denote this version by $\GMSOL$, for Guarded Monadic Second Order Logic.
We assume the reader is familiar with Monadic Second Order Logic $\MSOL$.
For the Monadic Second Order theory of graphs the reader is referred to the encyclopedic \cite{bk:CourcelleEngelfriet2011}.
$\CMSOL$ is the logic extending $\MSOL$ with modular counting quantifiers.

A {\em graph polynomial $P(G,\bar{X})$} is a {\em function $P$} which associates with
a graph $G$ a polynomial $P(G,\bar{X}) \in \Z[\bar{X}]$, and which is {\em invariant}
under graph isomorphisms (disregarding the labels).
Here $\bar{X} =(X_1, \ldots , X_s)$.

\medskip
Let us look at some examples:
\begin{enumerate}[(i)]
\itemsep=0.9pt
\item
Let $i_k(G)$ denote the number of independents sets $A \subseteq V(G)$ of size $k$.
\\
The {\em independence polynomial $In(G,X)$} is defined as
$$
In(G,X)= \sum_{k=0}^{n(G)} i_k(G) X^k
$$
\item
Let $c_k(G)$ denote the number of sets $A \subseteq V(G)$ which induce a clique of size $k$.
\\
The {\em clique polynomial $Cl(G,X)$} is defined as
$$
Cl(G;X) =\sum_{k=0}^{n(G)} c_k(G) X^k
$$
\item
Let $\lambda \in \N$. $\chi(G,\lambda)$ denotes the number of proper $\lambda$-colorings of $G$.
\\
By the well-known observation of G. Birkhoff (1912), $\chi(G,\lambda)$ is a polynomial in $\N[\lambda]$.
\\
We denote by $\chi(G,X)$ the extensions of $\chi(G,\lambda)$ to a polynomial over the complex numbers,
$\chi(G, X) \in \C[X]$.
$\chi(G,X)$ is called the {\em chromatic polynomial of $G$}.
\item
The {\em Tutte polynomial $T(G;X,Y)$} is defined as
\begin{gather}
T(G;X,Y) = \sum_{A \subset E(G)} (X-1)^{\kappa(A) - \kappa(G)} \cdot (Y-1)^{|A| + \kappa(A) -|V(G)|}
\label{tutte-def}
\end{gather}
where $\kappa(S)$ is the number of connected components of the spanning subgraph $G[S]=(V(G), S)$.
\item
The {\em matching polynomials} come in two versions:
Let $G$ be a graph and $m_k(G)$  be the number of matchings of size $k$ of $G$.
The {\em generating matching polynomial} $M(G:X)$ of $G$ is defined as
$$
M(G;X) = \sum_{k=0}^{\lfloor n/2 \rfloor} m_k(G) X^k
$$
and the {\em matching defect} aka {\em acyclic polynomial} $\alpha(G;X)$ is defined as
$$
\alpha(G;X) = \sum_{k=0}^{\lfloor n/2 \rfloor} (-1)^k m_k(G) X^{n-2k}
$$
The two are related by the equations
$$
\alpha(G;X) = X^n M(G;-X^{-2}) \text{  and  }
M(G;X)=  (-i)^n X^{n/2} \alpha(G;i X^{-1/2})
$$
There is also a bivariate version
$$
\bar{M}(G;X,Y) = \sum_{k=0}^{\lfloor n/2 \rfloor} (X)^k m_k(G) Y^{n-2k}
$$
where $\bar{M}(G;-1,Y)=\alpha(G;Y)$ and $\bar{M}(G;X,1)=M(G;X)$.
\item
In \cite{averbouch2008most,averbouch2010extension}
the authors introduce a most general edge elimination polynomial in three indeterminates $\xi(G;X,Y,Z)$.
\begin{gather}
\xi(G;X,Y,Z) = \sum_{(A \sqcup B) \subseteq E}
X^{\kappa(A \sqcup B) - c(B)}
\cdot
Y^{|A|+|B| - c(B)}
\cdot
Z^{c(B)}
\label{xi-def}
\end{gather}
Here $c((V,E))$ is the number of connected components of $(V,E)$
which have at least one edge.
Both the matching polynomials and the Tutte polynomial are substitution instances of $\xi(G;X,Y,Z)$.
In \cite{trinks2011covered,trinks2012proving}
other trivariate graph polynomials are discussed which are equivalent to $\xi(G;X,Y,Z)$:
Among them the subgraph counting polynomial $S(G,X,Y,Z)$ and the covered components polynomial $C(G;X,Y,Z)$.
All these trivariate graph polynomials are
substitution instances of each other.
The subgraph counting polynomial $S(G;X,Y,Z)$
is defined as\smallskip
$$
S(G;X,Y,Z) = \sum_{H=(W,F) \subseteq G}
X^{|W|}
Y^{\kappa(H)}
Z^{|F|}
$$
where $H$ ranges over all subgraphs of $G$.

\vspace{1mm}
The covered components polynomial $C(G;X,Y,Z)$ is defined as
$$
C(G;X,Y,Z) = \sum_{A\subseteq E}
X^{\kappa(([n],A))}
Y^{|A|}
Z^{c(([n],A))}
$$
\end{enumerate}
In \cite{trinks2011covered} it is shown how  $\xi(G;X,Y,Z)$ and $C(G;X,Y,Z)$ are related:
\begin{proposition}
\label{trinks}
$$
C(G;X,Y,Z) = \xi((G;X,Y, XYZ-XY)
\text{   and   }
\xi((G;X,Y,Z) = C(G;X,Y, \frac{Z}{XY} +1)
$$
\end{proposition}
A similar relation is given for $\xi(G;X,Y,Z)$ and $S(G;X,Y,Z)$ in \cite{trinks2012proving} which is a bit more complicated,
but no needed for this paper.

\begin{remark}
We note that both the Tutte polynomial and $\xi(G;,X,Y,Z)$ involve negative exponents,
whereas $S(G;X,Y,Z)$ and $C(G;X,Y,Z)$ do not.
This is the source of the difference between Theorem \ref{th:1} and Theorem \ref{th:4}.
\label{negative-exp}
\end{remark}

\subsection{The case of $G =K_n$}

Here we shall be concerned with computing a graph polynomial $P(G;\bar{X})$ in $k$ indeterminates
for the case where the graph $G$ is $K_n$, the complete graph on $n$ vertices.
We define, for
fixed non-negative integers $\bar{b} \in \N^k$
the sequence
$$P_{n}(\bar{b})= P(K_n,\bar{b})$$
for fixed values $\bar{b} \in \Z^k$.
Can we make some general statement about the sequence $P_n(\bar{b})$?
In some cases computing $P_n(\bar{b})$ is very easy, using trivial observations or simple recurrence relations.
However, the resulting graph polynomials may be unexpectedly complicated.
The Tutte polynomial and the
matching polynomial will illustrate this in the sequel.

\medskip
We compute $P_n(\bar{b})$ first for some straightforward cases:
\begin{enumerate}[(i)]
\itemsep=0.85pt
\item
{\em Independence polynomial}:
$$
i(K_n,k) =
\begin{cases}
0 & k \geq 2 \\
n & k=1 \\
1 & k=0
\end{cases}
$$
\eject

\noindent Hence,
$$
In(K_n, b) = i_1(K_n) \cdot b + c_0(K_n) = nb +1
$$
\item
{\em Clique polynomial}: $c_k(K_n) = {n \choose k}$, hence
$$
Cl(K_n,b) = \sum_{k=0}^n {n \choose k} b^k = (b+1)^n
$$
\item
{\em Chromatic polynomial}:
$\chi(K_n,b) = 0$ for $n\geq b+1$.
\end{enumerate}

\subsection{The Tutte polynomial}

The case of the Tutte polynomial is a bit more complicated, but of special interest.
We note that
\begin{enumerate}[(i)]
\itemsep=0.85pt
\item
$T(K_n;2,1)$
counts the number of forests on $n$ vertices.
\item
$T(K_n;1,1)$
counts the number of trees on $n$ vertices.
\item
$T(K_n;1,2)$
counts the number of connected graphs on $n$ vertices.
\end{enumerate}

The earliest computation of  $T_n(X,Y) = T(K_n,X,Y)$
can be found in \cite{tutte1954contribution} which already gives a recursive computation.
I. Gessel
\cite{gessel1979depth,gessel1995enumerative,gessel1996tutte}
and independently I. Pak \cite{pak-tutte}, proved:
$$
T(K_n; a,b) = \sum_{k=1}^n {(n-1) \choose (k-1)} \left(a + \sum_{i=1}^{k-1} b^i \right) T(K_{k-1};1,b) \cdot T(K_{n-k};a,b)
$$
I. Pak  in \cite{pak-tutte} also lists many other evaluations of $T(K_n; X,Y)$
with their combinatorial interpretations.

However,
in \cite{ar:BiggsDamerellSand72} it is noted that $T_n(X,Y)$
does not satisfy a linear recurrence which is independent of $n$.

\begin{proposition}[N.L. Biggs, R.M. Damerell and D.A. Sand]
There is no linear recurrence relation which computes the sequence $T_n(a,b)$ for fixed $a,b \in \Z$.
\end{proposition}

\subsection{The matching polynomial}

From \cite{godsil1981hermite} we know that for the defect matching polynomial (aka the acyclic polynomial)
$$
\alpha(K_n;X) = He_n(X)
$$
where $He_n(X)$ are the probabilist's Hermite polynomials for $n \in \N$.
The proof of this is due to C. Heilmann and E. Lieb
\cite{HeilmannLieb}.
From \cite[Equations 3.4 and 3.8]{carlitz1953congruence} one can derive\footnote{
Thanks to V. Rakita for checking this.
} that
the polynomials $He_n(X)$ satisfy the modular recurrence relation
$$
He_{n+m}(X) = He_n(X) \cdot  He_{m}(X) = He_n(X) \cdot X^m \left(\hspace*{-3mm}\mod{\mu}\right).
$$
and with $He_{0}(X) = 1$ and one gets  $He_m(X) =X^m$.

\begin{proposition}[Carlitz, 1953]
\label{pr:carlitz}
For every $a \in \N$ the sequence
$\alpha(K_n;a) = He_n(a) =a^n$ is ultimately period modulo $\mu$.
\end{proposition}
Our Theorem \ref{th:3} below shows the ultimate periodicity modulo $\mu$ for $M(K_n,a)$
without using the connection to the Hermite polynomials.

\subsection{Computing $P_n(\bar{b})$ modulo an integer $\mu$}

Let $\mu \in \N$.
We compute $P_n(\bar{b})$ modulo $\mu$ and observe:
\begin{enumerate}[(i)]
\itemsep=0.9pt
\item
For every $b \in \Z$ and $\mu \in \N$ the sequence
$In(K_n, b)= nb+1$ is {\em ultimately periodic}.
\item
For every $b \in \Z$ and $\mu \in \N$ the sequence
$Cl(K_n,b) =(b+1)^n$ is {\em ultimately periodic}.
\item
The sequence $\chi(K_n,b)$ is ultimately constant, hence  {\em ultimately periodic}.
\end{enumerate}
We shall see that this is the case
for a very large class of graph polynomials
subject to a definability condition in $\MSOL$.

Our main results are for the Tutte polynomial $T(G;X,Y)$, the bivariate matching polynomial $\bar{M}(G;X)$
and the trivariate edge elimination polynomial $\xi(G;X,Y,Z)$.
\begin{theorem}
\label{th:1}
For every $a,b, \mu \in\N^{+}$ with $a,b>1$,
$\gcd(a-1, \mu )=1$ and $\gcd(b-1, \mu )=1$, \\
the sequence $T(K_{n},a,b)$ is ultimately periodic modulo $\mu$.
\end{theorem}
Similarly, we also get for
the edge elimination polynomial $\xi(G;X,Y,Z)$,
the subgraph counting polynomial $S(G;X,Y,Z)$, and
the covered components polynomial $C(G;X,Y,Z)$:
\begin{theorem}
\label{th:4}
\begin{enumerate}[(i)]
\item
For every $a,b, c, \mu \in\N^+$
the sequences
$S(K_{n},a,b,c)$ and
$C(K_{n},a,b,c)$
are ultimately periodic modulo $\mu$ for every $\mu$.
\item
Assume  $ab$ divides $c$.
\\
Then $\xi(K_{n},a,b,c)$ is ultimately periodic modulo $\mu$ for every $\mu$.
\end{enumerate}
\end{theorem}
\begin{remark}
For $\xi(K_{n},a,b,c)$ one needs an additional condition like in Theorem \ref{th:1} due to the negative
exponents in the definition
\ref{xi-def}.
of $\xi(G;,X,Y,Z)$.
\end{remark}
\begin{theorem}
\label{th:3}
For every $a,b, \mu \in\N^{+}$
the sequence $\bar{M}(K_{n},a,b)$ is ultimately periodic modulo $\mu$.
\end{theorem}
Theorem \ref{th:3} can be generalized:
\begin{theorem}
\label{th:2}
If $P$ is a graph polynomial definable in $\GMSOL_{hgraph}$ without order,
then for every $\mu \in \N$ and every $a,b \in \Z$
the sequence $P(K_n, a,b)$ is ultimately periodic modulo $\mu$.
\end{theorem}
\begin{remark}
We note that Theorem \ref{th:1} is not a special case of Theorem \ref{th:2},
as the Tutte polynomial seems not to be $\GMSOL_{hgraph}$-definable without an order on the vertices.
\end{remark}

The three theorems only assert that the sequences are ultimately periodic modulo $\mu$, without
any indication of the length of the periodicity or the initial segment before the periodicity starts.

\subsection{Methods}

The proofs use tools from logic and combinatorics.
In particular
\begin{itemize}
\itemsep=0.9pt
\item
{\em Definability in Monadic Second Order Logic $\MSOL$};
\item
{\em Definability} in the extension $\CMSOL$ of $\MSOL$, where modular counting quantifiers are added;
\item
The {\em Specker-Blatter Theorem},
which gives a sufficient condition on when certain $\CMSOL$-definable density functions
are ultimately periodic if considered modulo $\mu$.
\end{itemize}

We assume that our readers are familiar with $\MSOL$ and first explain the {\em Specker-Blatter Theorem}.

\subsection{Related results}

The particular sequence
$T_n(1,b)=T(K_n;1,b)$
has been analyzed by
A. P. Mani and R.J. Stones,
\cite{ar:ManiStones2016},
for
$\mu = p^k$ where $p$ is an odd prime and $k \in \N$.

Let $\phi(n)$ be the {\em Euler totient function}
which is defined as the number of integers $a \in \{1, 2, . . . , m\} =[m]$ such that
$gcd(a, m) = 1$,  cf. \cite{bk:GrahamKnuthPatashnik94}.

\begin{proposition}[A.P. Mani and R.J. Stones]
Let $p$ be a prime, and let $k$ be a positive integer. For $b, n \in \Z$
such that $n \geq p^k$ and $b \neq 1 \mod{p}$, it holds that
$$
T_n(1,b) =\mod{p^k}
\begin{cases}
b^{\frac{\phi(p^k )}{2}} C_{n - \phi(p^k}(b) & \text{ if }  p \geq 3 \text{ and } n > p \\
b^{\frac{\phi(p)}{2}} - 1                   & \text{ if }  p \geq 3 \text{ and } n = p \\
1                                  & \text{ if }  p = n = 2\\
2                                  & \text{ if }  p = k = 2 \text{ and } n = 4\\
0                                  & \text{ otherwise}
\end{cases}
$$
\end{proposition}
This is much more informative than our Theorem \ref{th:1} for the case $a=1$ and $\mu = p^k$.
In \cite{ar:ManiStones2016} they formulate also a conjecture for $T(K_n;a,b)$ for general $a$, but still for $\mu = p^k$.

\begin{conjecture}[A.P. Mani and R.J. Stones]
Let $p$ be an odd prime and let $k$ be a positive integer.
\begin{enumerate}[(i)]
\itemsep=0.9pt
\item
If $n \geq p^k$, $a,b \in \Z$, and $b \not\equiv 1 \mod{p}$, then
modulo $p$ we have
$$
T(K_n;a,b) \equiv\mod{p}
\begin{cases}
b^{\frac{\phi(p)}{2} - 1} & \mbox{   if   } p\geq3, n=p  \\
& \mbox{ and } a=1 \mod{p}\\
b^{\frac{\phi(p^k)}{2}} \cdot T(K_{n - \phi(p^k)};a,b) & \mbox{otherwise}
\end{cases}
$$
\item
If $n \geq p^k$, $a,b \in \Z$, and $b \equiv 1 \mod{p}$, then
$$
T(K_n;a,b) \equiv \mod{p^k}
\begin{cases}
(n+a-1)^{p^k} \cdot T(K_{n -p^k};a,b) & n > p^k\\
(a-1)^{p^k-1} & n = p^k
\end{cases}
$$
\end{enumerate}
\end{conjecture}

\subsection{Outline of the paper}

In Section \ref{se:sb} we present the Specker-Blatter Theorem.
In Section \ref{se:matching} we prove Theorem \ref{th:3},
and in
In Section \ref{se:tutte} we prove Theorem \ref{th:1} and \ref{th:4}.
In Section \ref{se:msol} we prove Theorem \ref{th:2}.
Finally, in Section \ref{se:conclu} we present our conclusions and suggestions
for further research.

\section{The Specker-Blatter Theorem}
\label{se:sb}

The Specker-Blatter Theorem from 1980 is the first
application of Logic to Combinatorial Counting. At the time the theorem was hardly noticed,
mostly due to its unlucky placement for publication, \cite{pr:BlatterSpecker81,pr:BlatterSpecker84,ar:Specker88}.
An easily accessible place to find proofs and a survey of further developments is
\cite{ar:FischerKotekMakowsky11}.\vspace*{-2mm}

\subsection{Counting graphs: The density function}

A  {\em labeled graph $G$} with $n$ vertices will have $V(G) =\{0,1, \ldots , n-1\} =[n]$.
There are $gr(n)=2^{n \choose 2} = 2^{\frac{n(n-1)}{2}}$
many graphs with $n$ vertices.

A graph property is a class of finite graphs closed under graph isomorphisms.
For a graph property $\mathcal{P}$ denote by $\mathcal{P}^n$
the graphs with $n$ vertices in
$\mathcal{P}$, and by
$d_{\mathcal{P}}(n) = \mid \mathcal{P}^n \mid$,
the number of graphs $G$ with $V(G)=[n]$ which are in $\mathcal{P}$.
$d_{\mathcal{P}}(n)$ is called the {\em  density function} of $\mathcal{P}$.
The density function counts {\em labeled} graphs.
Let $G$ consist of the vertices $[n]$ and one single edge which is not a loop.
There is,$\,$up to isomorphisms (disregarding the labels) one such graph,$\,$but there are $n(n-1)/2$ such labeled$\,$graphs.

A graph property $\mathcal{P}$ is hereditary if it is closed under induced subgraphs.
$\mathcal{P}$ is monotone if it is closed under (not necessarily induced) subgraphs.
If
$\mathcal{P}$ is hereditary or monotone, the density function $d_{\mathcal{P}}(n)$  of
$\mathcal{P}$ is also called the
{\em  speed}  of
$\mathcal{P}$, since it is an {\em ultimately monotone increasing function},
cf. \cite{ar:BaloghBollobasWeinreich2000,ar:BaloghBollobasWeinreich2001,balogh2002measures}.
Studying the possible growth rate of the speed of a graph property was initiated in 1994
by E. Scheinerman and J. Zito in \cite{scheinerman1994size}.
E. Specker and C. Blatter already in 1980 studied under what conditions the density function of graph properties
satisfy recurrence relations.
The definition of the density function can be extended to relational structures of
any vocabulary $\tau$. However, we will consider only vocabularies $\tau$ without function symbols.
\begin{examples}
\begin{enumerate}[(i)]
\itemsep=0.9pt
\item
If $\mathcal{P} = Graphs$  consists of all  simple graphs,
$$
d_{Graphs}(n) = 2^{n \choose 2}
$$
In the unlabeled case the function is rather complicated.
\item
If $\mathcal{P} = LinOrd$ consists of all linear orders.\vspace*{-1mm}
$$
d_{LinOrd}(n) =  n!
$$
In the unlabeled case we have  the constant function
with value $1$.
\item
If $\mathcal{P} = SqGrids$ consists of all square grids,\vspace*{-1mm}
$$
d_{SqGrids}(n) =
\begin{cases}
\frac{n!}{4} & \mbox{ if } n = m^2 \\
0 & \mbox{ else }
\end{cases}
$$
In the unlabeled case we have $1$ instead of $n!$ in the above expression.
\end{enumerate}
\end{examples}

\vfil\eject

For some graph properties $\mathcal{P}$
the density functions satisfies a linear recurrence relation over $\Z$.
However, this is not always the case.

\begin{lemma}[Folklore]
\label{le:rec-growth}
Let $f: \mathbb{Z} \rightarrow \mathbb{Z}$ a function which satisfies
a linear recurrence relation
$$f(n+1) = \sum_{i=0}^k a_i f(n-i)$$
over $\mathbb{Z}$.
Then
there is a constant $c \in \mathbb{Z}$ such that $f(n) \leq 2^{cn}$.
\end{lemma}

\begin{corollary}
For $\mathcal{C} \in \{ Graphs, LinOrd, SqGrids\}$,
$d_{\mathcal{C}}(n)$ does not satisfy a linear recurrence over
$\mathbb{Z}$.
\end{corollary}

\subsection{Modular counting}

Let $\mu \in \N$.
\begin{obs}
For every $\mu \in \mathbb{N}$ and for large enough $n$ we have
$n! =0 \pmod{\mu}$
\end{obs}
Hence, for $n \geq N(\mu)$ we have
$$
d_{LinOrd}(n+1) = d_{LinOrd}(n) \pmod{\mu}\vspace*{-2mm}
$$
and\vspace*{-1mm}
$$
d_{SqGrid}(n+1) = d_{SqGrid}(n) \pmod{\mu}\vspace*{1mm}
$$
We say that a function $f(n)$ satisfies a {\em trivial modular recurrence}
if for every $\mu$ there exists $N_{\mu}$ such that
if $n>N_{\mu}$ then $f(n)\equiv 0 \pmod{\mu}$.
Clearly,
the two examples above satisfy
trivial modular recurrences.

\begin{obs}
$f(n)$ satisfies a trivial modular recurrence
iff
there exist functions $g(n), h(n)$ with $g(n)$ tending to infinity
such that $f(n) = g(n)! \cdot h(n)$.
\end{obs}
In other words,
trivial modular recurrences are always caused by some factor which is a factorial.

\medskip
It is sometimes more intuitive to say of integer sequences with values in $[m]$ that they
are ultimately periodic rather than to talk about modular linear recurrence.

\begin{proposition}[Folklore]
Let $a_n$ be an integer sequence.
\\
$a_n$ satisfies a linear recurrence relation modulo $\mu$
iff
$a_n$ is ultimately periodic modulo $\mu$.
\end{proposition}

The following are two instructive examples.
First, we look at the class of all graphs.
\begin{example}
\label{ex:1}
The density function for all graphs is given by
$$
d_{Graphs}(n+1) =
2^{n+1 \choose 2} =
2^{n \choose 2} \cdot 2^n.
$$
Therefore\vspace*{-2mm}
$$
d_{Graphs}(n+m) =
d_{Graphs}(n)  \cdot \prod_{i=0}^{m-1} 2^{n+i} =
d_{Graphs}(n)  \cdot 2^{nm} \cdot \prod_{i=0}^{m-1} 2^{i}
$$
As $a^{p-1}  =a \pmod{p}$ (Fermat's Little Theorem) we get
with $a=2^n$ and $m=p$ a prime
$$
d_{Graphs}(n+p) =
d_{Graphs}(n)  \cdot \prod_{i=0}^{p-1} 2^{i} \pmod{p}
$$
This is a non-trivial recurrence for $\mu=p$ a prime.
It is also different for distinct primes $p$ and $p'$,
In other words, the existence of the modular recurrence is non-uniform in $p$.
\end{example}

The second example is the class of graphs
$EQ_2CLIQUE$
which consists of the graphs which are the disjoint unions of two equal-sized cliques.

\begin{example}
\label{ex:2}
For density function $d_{EQ_2CLIQUE}(n)$ we have
$$
d_{EQ_2CLIQUE}(n) = b_2(n)=
\begin{cases}
\frac{1}{2} {2m \choose m} & \mbox{ for } n=2m
\\
0 & \mbox{ else }
\end{cases}
$$
The factor $\frac{1}{2}$ is there because we cannot distinguish
the choice of the first clique from the choice of its complement.
\end{example}

The function $b_2(n)$ was studied by
F. E. A. Lucas (1842-1891)
in 1878, but not published at the time.
It was found in his notes in the National Archive of France.
A proof may be reconstructed from the hints in \cite[Exercise 5.61]{bk:GrahamKnuthPatashnik94}.
\begin{proposition}[Lucas, 1878]
For every $n$ which is not a power of $2$, we have
$b_2(n)\equiv 0 \pmod{2}$, and for every $n$ which is a power
of $2$ we have $b_2(n)\equiv 1 \pmod{2}$.
\\
In particular, $b_2(n)$ is not ultimately periodic modulo $2$.
\end{proposition}
We conclude that
$d_{EQ_2CLIQUE}(n)$ is not ultimately periodic modulo $2$.

Finally, here is an example, were the precise counting is known. It is originally due to
Redfield, \cite{ar:Redfield27} and was rediscovered by R.C. Read and G. Polya, \cite{ar:Read59,ar:Read60,bk:PolyaRead}.
\begin{example}
\label{ex:3}
Let $\mathcal{R}_d$ be the class of regular graphs of degree $d$
and $d_{{\mathcal R}_d}(n)$ its density function.
\end{example}
It is not at all obvious that $d_{{\mathcal R}_d}(n)$ is ultimately periodic modulo $\mu$, but it follows
from the Specker-Blatter Theorem \ref{th:sb} below that, indeed, it is.
\begin{proposition}[J.H. Redfield, 1927]
For $d=3$ we have
$d_{{\mathcal R}_3}(2n+1) = 0$ and
$$
d_{{\mathcal R}_3}(2n) =\frac{(2n)!}{6^n}?\sum_{j,k}
\frac{(-1)^j(6k-2j)!6^j}{(3k-j)!(2k-j)!(n-k)!}48^k\sum_i\frac{(-1)^i j!}{(j-2i)! i!}
$$
\end{proposition}
An accessible proof can be found in
\cite[page 187]{bk:HararyPalmer}.

\subsection{$\MSOL$-definable graph properties}

We now look at the density functions of graph properties $\mathcal{P}$
which are definable in $\MSOL$.

\begin{theorem}[Specker-Blatter Theorem]
\label{th:sb}
Let $\mathcal{P}$
be a graph property which is $MSOL$-definable
and let
$d_{\mathcal{P}}(n)$ be its density function.
\begin{itemize}
\itemsep=0.95pt
\item
$d_{\mathcal{P}}(n)$ satisfies modular recurrence relations for
each $\mu \in \N$, hence it is ultimately periodic modulo $\mu$.
\item
This remains true for vocabularies $\tau$ with several binary edge relations and unary predicates
on the vertices.
\end{itemize}
\end{theorem}

The Specker-Blatter Theorem does not hold
if one allows quaternary relations in $\tau$,
\cite{ar:Fischer02,ar:Specker2005}.
\begin{theorem}[E. Fischer, 2003]
Let $\tau_0$ consist of one quaternary relation.
There is a class of
$FOL(\tau_0)$-definable  $\tau_0$-structures $\mathcal{F}$
such that
$d_{\mathcal{F}}$ is not ultimately periodic modulo $2$.
\end{theorem}
The proof consists of a very clever encoding of
$EQ_2CLIQUE$ using the quaternary relation.

The restriction to binary relations is not needed if the graph property $\mathcal{P}$ contains only graphs of bounded degree,
see \cite{pr:FischerMakowsky03}. The same works for classes of $\tau$-structures of bounded degree, where
the degree is defined via the Gaifman graph of the structures.

\medskip
The Gaifman graph of a $\tau$-structure $\aA$ is the (undirected, loop-free) graph $G_{\aA}$ with
vertex set $A$, the universe of $\aA$ and an edge between two distinct vertices $a, b in A$ iff there exists an $R  \in \tau$
and a tuple $(a_1, . . , a_r) \in R^A$ such that $a, b \in \{a_1, . . , a_r\}$, cf. \cite{bk:Libkin2004}.

\section{The matching polynomial of the complete graph}
\label{se:matching}

Here we prove
\begin{quote}
{\bf Theorem \ref{th:3}}:
For all $a,b, \mu \in \N$.
$\bar{M}(K_n;a,b)$
is ultimately periodic modulo $\mu$.
\\
In particular $M(K_n;a)$
is ultimately periodic modulo $\mu$.
\end{quote}
\begin{proof}
First we prove it for $M(K_n;a)$ and note first that $M(K_n;a)$
is of the form
$$
M(K_n;a) =
P_{\psi}(K_n,a) = \sum_{F \subseteq E(G): \psi(E,F)} a^{|F|}
$$
where $\psi(E,F)$ says that $E$ is the edge relation of $K_n$ and $F \subseteq E$ is a set of independent edges of $K_n$.
We interpret $a^{F}$ as the set of functions $f: F \rightarrow [a]$.
Each $f$ induces a  partition of $F$ with
$$F_f(i) = \{ e \in F: f(e) =i \}.$$
We have
$$
a^{|F|} = |\{ f: F \rightarrow [a] \}|.
$$
Next
we look at the density function of
$$
g_{a}(n) = |\{ U_1, \ldots U_a \subseteq [n]^2, F \subseteq [n]^2 : \phi_1(\bar{U},F), \psi(E,F) \} |
$$
where
$\phi_1(\bar{U})$ says: $U_1, \ldots, U_{a}$ partition $F$ and
$\psi(E,F)$ says that $E$ is the edge relation of $K_n$ and $F \subseteq E$ is a set of independent edges of $K_n$.
$g_{a}(n)$ encodes the computation of $M(K_n;a)$.
\begin{claim}
$$
g_{a}(n) = P_{\psi}(K_n,a) = \sum_{F \subseteq E(G): \psi(E,F)} a^{|F|} =M(K_n;a)
$$
\end{claim}
All the relation symbols of $\phi_1$ and $\psi$
are binary, therefore we can apply the Specker-Blatter Theorem and conclude that
$g_{a}(n)$ is ultimately periodic modulo every $\mu \in \N$, and so is $M(K_n;a)$.

\medskip
The proof for $\bar{M}(K_n;a,b)$ is similar.
$$
\bar{M}(K_n;a,b) =
\sum_{k=0}^{\lfloor n/2 \rfloor} (X)^k m_k(G) Y^{n-2k} =
\sum_{F \subseteq E(G): \psi(E,F)}
(a)^{|F|} m_k(G) b^{|dF|}
$$
where $dF$ is the set of vertices in $[n]$ not covered by $F$.
If $F$ has $k$ edges, then $dF$ has $n-2k$ vertices.
\end{proof}

\vspace*{-3mm}
\section{The Tutte polynomial of a complete graph}
\label{se:tutte}

Now we prove
\begin{quote}
{\bf Theorem \ref{th:1}}
\\
For every $a,b, \mu \in\N^{+}$ with $a,b>1$,
$\gcd(a-1, \mu )=1$ and $\gcd(b-1, \mu )=1$,
$T_n(a,b)=T(K_{n},a,b)$ is ultimately periodic modulo $\mu$.
\end{quote}

\begin{proof}
We first rewrite the Tutte polynomial as
$$
T(G;X,Y) =
\frac{1}{(Y-1)^{|V(G)|} \cdot (X-1)^{\kappa(G))}}
\sum_{A \subseteq E(G)} (X-1)^{\kappa(A)} \cdot (Y-1)^{\kappa(A) + |A|}
$$
and put $G =K_n$. Now $|V(K_n)|=n$ and $|E(K_n)|={n \choose  2}$, and  $\kappa(K_n)=1$.
$$
T(K_n;X,Y) =
\frac{1}{(Y-1)^{n} \cdot (X-1)}
\sum_{A \subseteq E(G)} (X-1)^{\kappa(A)} \cdot (Y-1)^{\kappa(A) + |A|}.
$$
\eject
Like in the proof Theorem \ref{th:3}
we interpret
$a^{A}$ as the set of functions $f: A \rightarrow [a]$.
Each $f$ induces a  partition of $A$ with
$$A_f(i) = \{ e \in A: f(e) =i \}.$$
We have
$$
a^{|A|} = |\{ f: A \rightarrow [a] \}|.
$$
In the case of
$a^{\kappa(A)}$ the sets $A_f(i)$ have to be $A$-closed, in order to partition the connected components
of the spanning induced by the set of edges of $A$.
Next look at the function
$$
f_{a,b}(n)= (b-1)^{n} \cdot (a-1) \cdot T(K_n, a,b)
=
\sum_{A \subseteq E(G)} (a-1)^{\kappa(A)} \cdot (b-1)^{\kappa(A)} \cdot (b-1)^{|A|}
$$
as a  $\MSOL$-definable density function:
\begin{gather}
f_{a,b}(n) = |\{\bar{U}, \bar{R}, \bar{S} \subseteq [n]:
\phi_1(\bar{U},A), \phi_2(\bar{R}, A), \phi_3(\bar{S}), A \subseteq E \} |
\notag
\end{gather}
where for
$\bar{U} = (U_1, \ldots, U_{a-1})$ and
$\bar{R} = (R_1, \ldots, R_{b-1})$ are unary relations and
$\bar{S} = (S_1, \ldots, S_{b-1})$ are binary, and
\begin{enumerate}[(i)]
\item each $U_1, \ldots, U_{a-1}, R_1, \ldots, R_{b-1}, S_1, \ldots, S_{b-1} \subseteq [n]$;
\item  $\phi_1(\bar{U}, A)$ says: $U_1, \ldots, U_{a-1}$
partitions the connected components of $G[A]$,
and each $U_i$ is a disjoint union of connected components of the graph $G= ([n], A)$;
\item
$\phi_2(\bar{R}, A)$ says: $R_1, \ldots, R_{b-1}$ also
partitions the connected components of $G = ( [n], A)$,
\item and
$\phi_3(\bar{S})$ says: $S_1, \ldots, S_{b-1}$ partitions $A \subseteq [n]^2$.
\end{enumerate}
All the formulas
$\phi_1(\bar{U}, A), \phi_2(\bar{R}, A), \phi_3(\bar{S})$
are in $\MSOL$
and
contain only unary and binary relation symbols.
It follows by Theorem \ref{th:sb} that $f_{a,b}(n)$ is
ultimately periodic modulo every $\mu \in \N$.
\begin{claim}
$$
f_{a,b}(n) = \sum_{A \subseteq E(G)} (a-1)^{\kappa(A))} \cdot (b-1)^{\kappa(A)} \cdot (b-1)^{|A|}
$$
\end{claim}

Now we need a lemma.
\begin{lemma}
\label{lemma}
Let $d_1(n), d_2(n)$ be integer functions.
\begin{enumerate}[(i)]
\item
Let $c , \mu \in \N^+ $.
\\
Assume $c \cdot d_1(n)$ is ultimately periodic modulo $\mu$ and $\gcd(c, \mu)=1$.
\\
Then $d_1(n)$ is ultimately periodic modulo $\mu$.
\item
Let $t, \mu \in \N^+$ with $\gcd(t,\mu)=1$.
\\
Assume that $t^n \cdot d_2(n)$ is ultimately periodic modulo $\mu$.
\\
Then $d_2(n)$ is ultimately periodic modulo $\mu$.
\end{enumerate}
\end{lemma}
\begin{proof}
(i) is left to the reader.
\\
(ii):
Since $t$ and $\mu$ are relatively prime, $t$ has a multiplicative
inverse $g$ modulo $\mu$: $g \cdot t \equiv 1 (\mod{\mu})$.
The product of two ultimately periodic sequence is ultimately
periodic, hence this is true for the product of $g$ and  $t^n \cdot d_2(n)$.
But we have
$$g^n \cdot t^n \cdot d_2(n) \equiv d_2(n)(\mod{\mu}),$$ hence $d_2(n)$
is ultimately periodic.
\end{proof}
To complete the proof of Theorem \ref{th:1} we note that $\kappa(E) =1$  and $|V|=n$,
and we put:
$$
d_1(n) = T(K_n, a,b) \mbox{  and  } d_2(n) = (b-1) \cdot T(K_n, a,b)
$$
We have $f_{a,b}(n)= (a-1)^n \cdot d_2(n)$ is ultimately periodic modulo $\mu$
for $\gcd((a-1),\mu)=1$.
By the Lemma  \ref{lemma}(ii) we have that $d_2(n)$ ultimately periodic modulo $\mu$
for $\gcd((a-1),\mu)=1$.
By the Lemma  \ref{lemma}(i) we have that $d_1(n)$ ultimately periodic modulo $\mu$
for $\gcd((b-1),\mu)=1$.
\end{proof}

The proofs of Theorem \ref{th:4}  for $S(G;a,b,c)$ and $C(G;a,b,c)$ are similar as for  $f_{a,b}(n)$.
\\
For $\xi(G;a,b,c)$ we note that $\xi((G;X,Y,Z) = C(G;X,Y, \frac{Z}{XY} +1)$ from Proposition \ref{trinks}.
\\
To be able to use now $C(G;a,b, \frac{c}{ab} +1)$ one requires that $ab$ divides $c$ so that
$\frac{c}{ab}$ is a non-negative integer.

\section{Theorem \ref{th:2} and its limitations}
\label{se:msol}

The logic $\CMSOL$ is the extension of $\MSOL$ by adding {\em modular counting quantifiers}.
Let $\phi(x)$ be a formula with free variable $x$.
A modular counting quantifier $C_{\mu, k}x \phi(x)$
says that there are, modulo $\mu$, $k$-many elements satisfying $\phi(x)$.
Using Ehrenfeucht-Fra\"iss\'e games for $\MSOL$ one can show that
$C_{\mu, k}x \phi(x)$ is not expressible in $\MSOL$. However,
for $\phi(x)$ a formula in $\SOL$ the formula
$C_{\mu, k}x \phi(x)$ it is expressible in $\SOL$.
We say that there is an equivalence relation $E$
on the set defined by $\phi(x)$ which has exactly one equivalence class of size $k$ and all the other
non-empty equivalence classes are of size $\mu$. There are other ways of defining $C_{\mu, k}x \phi(x)$
in $\SOL$. Instead of the existence of one binary relation we can also assert that there are $\mu$ disjoint subsets
of equal size
of the set defined by $\phi(x)$
and the complement of their union has size $k$.
But then we need the existence of a binary relation which expresses
that the unary predicates are of equal size.

\medskip
$\CMSOL$ is like $\MSOL$ but modular counting quantifiers are allowed.
The syntax of $\CMSOL$ is obtained from the syntax of $\MSOL$ by allowing
also quantification with $C_{\mu, k}x$. The meaning function of $\MSOL$ is then
naturally extended to $\CMSOL$.
$\CGMSOL$ is defined analogously by adding modular counting quantifiers to $\GMSOL$.

It was shown in \cite{pr:FischerMakowsky03,ar:FischerKotekMakowsky11}
that
the Specker-Blatter Theorem also holds for $\CMSOL$ for relational structures with
relations of arity at most two, and for $\CGMSOL$ graphs.

\subsection{$\CMSOL$-definable graph polynomials}

We now discuss how Theorems \ref{th:1} and \ref{th:3} can be extended.

\medskip
A univariate graph polynomial of the form
$$
P_{\phi}(G;X) =\sum_{A \subseteq V(G): \phi(E,A)} X^{|A|}
$$
is $\CMSOL$ definable if $\phi(A,E)$
is a $\CMSOL$-formula in the language of graphs with an additional predicate for $A$.

\medskip
A univariate graph polynomial of the form
$$
P_{\psi}(G;X) =\sum_{F \subseteq E(G): \psi(E,F)} X^{|F|}
$$
is $\CGMSOL$ definable if
$\psi(F,E)$
is a $\CGMSOL$-formula in the language of graphs with an additional predicate for $A$
or for $F \subseteq E$.
The independence polynomial and the clique polynomial
are of the first form.
The matching polynomials $\alpha$ and $M$ are of the second form.

\medskip
Analyzing the proof of Theorem \ref{th:3} immediate gives Theorem \ref{th:2}:
\begin{quote}
If $P$ is a graph polynomial definable in $\CGMSOL_{graph}$ without order,
then for every $\mu \in \N$ and every $a,b \in \Z$\\
the sequence $P(K_n, a,b)$ is ultimately periodic modulo $\mu$.
\end{quote}
The proof of Theorem \ref{th:1} also works for the graph polynomials listed in Theorem \ref{th:4}.
Unfortunately, we have not found interesting applications of Theorem \ref{th:2}.
In many cases $P(K_n, a,b)$ can be shown directly to
satisfy linear recurrence relations over $\Z$.
\begin{examples}
\begin{enumerate}[(i)]
\itemsep=0.9pt
\item
If $\mathcal{P}$ is a graph property which does not contain any complete graph and is
definable in $\CMSOL$ for graphs
and $\phi_{\cP}(A,E)$ says that $A$ induces a graph $G[A] \in \cP$
then
$$
P_{\phi_{\cP}}(G;X) =\sum_{A \subseteq V(G): \phi_{\cP}(E,A)} X^{|A|}
$$
trivializes for $P_{\phi_{\cP}}(K_n;X) $.
\item
If instead, $\phi^{'}_{\cP}(A,E')$ says that $(A, E')$ is a subgraph  $\cP$
$P_{\phi^{'}_{\cP}}(K_n;X)$ may be non-trivial.
\item
The domination polynomial $D(G;X)$ is obtained by taking
$\phi_{dom}(A,E)$ which says that $A$ is a dominating set for $E$.
In \cite{alikhani2009introduction,oboudi2016roots}
it is shown that $D(K_{n+1};X) = D(K_n;X)(X+1) +X$.
For $X =a$ this is a linear recurrence relation, hence Theorem \ref{th:2}
gives nothing new.
\item
The univariate interlace polynomial $q(G;X)$ from \cite{arratia2004interlace} is $\CMSOL$-definable, as shown in
\cite{courcelle2008multivariate}. But  $q(K_n;X) = 2^{n-1}X$, hence $q(K_{n+1};X) = 2X q(K_n;X)$.
Again Theorem \ref{th:2} gives nothing new.
\end{enumerate}
\end{examples}

\begin{problem}
Find  more graph polynomials $P(G; \bar{a})$  where for $\bar{a} \in \N^k, \mu \in \N$ the
sequence $P(K_n; \bar{a})$ is not obviously ultimately periodic modulo $\mu$.
\end{problem}

\begin{problem}
Find interesting cases of graph polynomials  where Theorem \ref{th:2} gives a non-trivial result.
\end{problem}

We have seen that Theorem \ref{th:2} also holds for multivariate polynomials.
The Tutte polynomial is not of this form, because of the term $(X-1)^{\kappa(S)}$.
It is of this form in the language of ordered graphs.
However, the Specker-Blatter Theorem formulated for ordered structures trivializes, because
the ordering adds a factor of the form $n!$, hence the sequence satisfies a trivial modular recurrence.

\subsection{Why complete graphs?}

Complete  labeled graphs on $n$ vertices are definable in $\FOL$ in the empty vocabulary
and are unique, not only up to isomorphisms.
The same is true about the empty graph.

\begin{proposition}
Assume $\phi(\bar{x},\bar{y})$ is a formula of $\CMSOL$ over the empty vocabulary which
defines a unique edge relation  $E_{\phi}$ on $k$-tuples of a set $[n]$.
Then the graph $G = ([n]^k,  E_{\phi})$ is either a complete graph or the empty graph.
\end{proposition}
The way we used the Specker-Blatter Theorem  in the proof of Theorem \ref{th:2} did require that the function
$$
g_{a}(n) = |\{ U_1, \ldots U_a \subseteq [n], F \subseteq [n]^2 : \phi_1(\bar{U},F), \psi(F) \} |
$$
evaluates the matching polynomial at $a$ for $G = K_n$.
If instead of $G = K_n$ we use some $G = ([n]^k,  E_{\phi})$ with $k=1$ we can modify $g_{a}(n)$
and still apply the Specker-Blatter Theorem to it, but the modified version $g_{a, \phi}(n)$ will
not evaluate the matching polynomial anymore.
Let $sp(\phi, n)$ be the number of labeled graphs of the form
$G = ([n],  E_{\phi})$.
Even if all the graphs
$G_n = ([n],  E_{\phi})$ are isomorphic, $g_{a, \phi}(n)$
computes
$$
g_{a, \phi}(n)=
M(G_n;a) \cdot sp(\phi, n)
$$
which will be ultimately periodic modulo $\mu$. However, this does not suffice to conclude
that
$M(G_n;a)$ is ultimately periodic modulo $\mu$.

\medskip
{\bf Facit:} Our proofs of Theorems
\ref{th:1},
\ref{th:2} and
\ref{th:3}
only work for $G=K_n$ or $G=\bar{K_n}$.

\section{Conclusions}
\label{se:conclu}

Inspired by the paper of A.P. Mani and R.J. Stones \cite{ar:ManiStones2016}
we have examined modular recurrences for the matching and the Tutte polynomial
of a complete graph. We have noted that the existence of modular recurrence relation
follows from the Specker-Blatter Theorem, without explicitly describing the
exact modular recurrence. The conjectures of A.P. Mani and R.J. Stones
are more ambitious, as they give a precise statement about how these modular recurrences
look in the case of the Tutte polynomial.
We also noted that our approach via the Specker-Blatter Theorem \ref{th:sb}
works for other graph polynomials definable in variants of $\MSOL$.
Most strikingly, it works for the trivariate edge elimination polynomial $\xi(G;X,Y,Z)$,
which is the most general edge elimination polynomial, and generalizes both the matching polynomials
and the Tutte polynomial.

Our result suggests that the use of the Specker-Batter Theorem to
$\CMSOL$-definable ($\CGMSOL$-definable) graph polynomials should be further
investigated.
In \cite{fischer2003specker,ar:FischerKotekMakowsky11} a logic-free version of the Specker-Batter Theorem is
discussed, where instead of the graph property $\cP$  being definable in $\CMSOL$
one only requires that $\cP$ has finite Specker rank.
Here we note that there uncountably many graph properties of finite Specker rank,
but there only countably many $\CMSOL$-definable graph properties.
$\CMSOL$-definable graph properties have always finite Specker rank, but the upper bound on the Specker rank
which follows from $\CMSOL$-definability is very often exponentially bigger than the true Specker rank.
On the other hand the parameters which characterize the ultimate periodicity depend only polynomially
on the Specker rank, hence
the structure of the modular recurrences can be more precisely described
using the Specker rank.

\subsection*{Acknowledgments}

The material of this paper was first presented at the
Colloquium Logicum, September 10-12, 2016 held in Hamburg, Germany.
An updated version was given at the Colloquium in honor of B. Trakhtenbrot's centenary,
held in Tel Aviv (and via ZOOM) in October 2021.
The second author also presented a survey of the Specker-Blatter Theorem and its
various extensions at the Colloquium held at ETHZ in honor of E. Specker's centenary in February 2020
in Zurich, Switzerland.
Both E. Specker and C. Blatter were my (the second author's) teachers in 1967-1973. They both shaped my
mathematical and logical views in a lasting way. Boris Abramovich did the same after we
both started to live in Israel after 1980. He gave his unrelenting moral support for my own work.
Thanks also to the unknown referee for his careful reading and his resulting suggestions, which helped improving the
paper. Thanks also to Peter Tittmann for useful comments.


\begin{thebibliography}{10}
\providecommand{\url}[1]{\texttt{#1}}
\providecommand{\urlprefix}{URL }
\expandafter\ifx\csname urlstyle\endcsname\relax
  \providecommand{\doi}[1]{doi:\discretionary{}{}{}#1}\else
  \providecommand{\doi}{doi:\discretionary{}{}{}\begingroup
  \urlstyle{rm}\Url}\fi
\providecommand{\eprint}[2][]{\url{#2}}

\bibitem{trakhtenbrot1962finite}
Trakhtenbrot BA.
\newblock Finite automata and monadic second order logic.
\newblock \emph{Siberian Math. J}, 1962.
\newblock \textbf{3}:101--131.
\newblock English translation in: AMS Transl. 59 (1966), 23--55.

\bibitem{ar:FischerMakowsky08}
Fischer E, Makowsky J.
\newblock Linear Recurrence Relations for Graph Polynomials.
\newblock In: Avron A, Dershowitz N, Rabinowitz A (eds.), Boris (Boaz) A.
  Trakhtenbrot on the occasion of his 85th birthday, volume 4800 of
  \emph{LNCS}. Springer, 2008 pp. 266--279.
  doi:10.1007/978-3-540-78127-1\_15.

\bibitem{pr:BlatterSpecker84}
Blatter C, Specker E.
\newblock Recurrence relations for the number of labeled structures on a finite  set.
\newblock In: B{\"o}rger E, Hasenjaeger G, R{\"o}dding D (eds.), In Logic and
  Machines: Decision Problems and Complexity, volume 171 of \emph{Lecture Notes
  in Computer Science}. Springer, 1984 pp. 43--61.

\bibitem{bk:CourcelleEngelfriet2011}
Courcelle B, Engelfriet J.
\newblock Graph Structure and Monadic Second-order Logic, a Language Theoretic
  Approach.
\newblock Cambridge University Press, 2012.
doi:10.1017/CBO9780511977619.

\bibitem{averbouch2008most}
Averbouch I, Godlin B, Makowsky JA.
\newblock A most general edge elimination polynomial.
\newblock In: International Workshop on Graph-Theoretic Concepts in Computer
  Science. Springer, 2008 pp. 31--42.
  doi:10.1007/978-3-540-92248-3\_4.

\bibitem{averbouch2010extension}
Averbouch I, Godlin B, Makowsky JA.
\newblock An extension of the bivariate chromatic polynomial.
\newblock \emph{European Journal of Combinatorics}, 2010.
\newblock \textbf{31}(1):1--17.   doi:10.1016/j.ejc.2009.05.006.

\bibitem{trinks2011covered}
Trinks M.
\newblock The covered components polynomial: A new representation of the edge
  elimination polynomial.
\newblock \emph{arXiv preprint arXiv:1103.2218}, 2011.

\bibitem{trinks2012proving}
Trinks M.
\newblock Proving properties of the edge elimination polynomial using
  equivalent graph polynomials.
\newblock \emph{arXiv preprint arXiv:1205.2205}, 2012.

\bibitem{tutte1954contribution}
Tutte WT.
\newblock A contribution to the theory of chromatic polynomials.
\newblock \emph{Canadian journal of mathematics}, 1954.
\newblock \textbf{6}:80--91.

\bibitem{gessel1979depth}
Gessel IM, Wang DL.
\newblock Depth-first search as a combinatorial correspondence.
\newblock \emph{J. Comb. Theory, Ser. A}, 1979.
\newblock \textbf{26}(3):308--313.
doi:10.1016/0097-3165(79)90108-0.

\bibitem{gessel1995enumerative}
Gessel IM.
\newblock Enumerative applications of a decomposition for graphs and digraphs.
\newblock \emph{Discrete mathematics}, 1995.
\newblock \textbf{139}(1-3):257--271.
doi:10.1016/0012-365X(94)00135-6.

\bibitem{gessel1996tutte}
Gessel IM, Sagan BE.
\newblock The Tutte Polynomial of a Graph, Depth-first Search.
\newblock \emph{The Electronic Journal of Combinatorics}, 1996.
\newblock pp. R9--R9.   doi:10.37236/1267.

\bibitem{pak-tutte}
Pak IM.
\newblock Computation of Tutte polynomials for complete graphs.
URL \texttt{https://www.math.ucla.edu/} \texttt{pak/papers/Pak\_Computation\_Tutte\_polynomial\_complete\_graphs.pdf.}

\bibitem{ar:BiggsDamerellSand72}
Biggs N, Damerell R, Sand D.
\newblock Recursive families of graphs.
\newblock \emph{J. Combin. Theory Ser. B}, 1972.
\newblock \textbf{12}(2):123--131.
doi:10.1016/0095-8956(72)90016-0.

\bibitem{godsil1981hermite}
Godsil CD.
\newblock Hermite polynomials and a duality relation for matchings polynomials.
\newblock \emph{Combinatorica}, 1981.
\newblock \textbf{1}(3):257--262.
doi:10.1007/BF02579331.

\bibitem{HeilmannLieb}
Heilmann C, Lieb E.
\newblock Theory of monomer-dymer systems.
\newblock \emph{Comm. Math. Phys}, 1972.
\newblock \textbf{25}:190--232.
doi:10.1007/BF01877590.

\bibitem{carlitz1953congruence}
Carlitz L.
\newblock Congruence properties of the polynomials of Hermite, Laguerre and
  Legendre.
\newblock \emph{Mathematische Zeitschrift}, 1953.
\newblock \textbf{59}(1):474--483.

\bibitem{ar:ManiStones2016}
Mani A, Stones R.
\newblock The number of labeled connected graphs modulo prime powers.
\newblock \emph{SIAM. J. Discrete Math.}, 2016.
\newblock \textbf{30.2}:1046--1057.
doi:10.1137/15M1024615.

\bibitem{bk:GrahamKnuthPatashnik94}
Graham R, Knuth D, Patashnik O.
\newblock Concrete Mathematics.
\newblock Addison-Wesley, 2 edition, 1994.

\bibitem{pr:BlatterSpecker81}
Blatter C, Specker E.
\newblock Le nombre de structures finies d'une th'eorie \`a charact\`ere fin.
\newblock \emph{Sciences Math\'ematiques, Fonds Nationale de la recherche
  Scientifique, Bruxelles}, 1981.
\newblock pp. 41--44.

\bibitem{ar:Specker88}
Specker E.
\newblock Application of Logic and Combinatorics to Enumeration Problems.
\newblock In: B{\"o}rger E (ed.), Trends in Theoretical Computer Science.
  Computer Science Press, 1988 pp. 141--169.
\newblock Reprinted in: Ernst Specker, Selecta, Birkh\"auser 1990, pp. 324-350.
doi:10.1007/978-3-0348-9259-9\_29.

\bibitem{ar:FischerKotekMakowsky11}
Fischer E, Kotek T, Makowsky J.
\newblock Application of Logic to Combinatorial Sequences and Their Recurrence  Relations.
\newblock In: Grohe M, Makowsky J (eds.), Model Theoretic Methods in Finite
  Combinatorics, volume 558 of \emph{Contemporary Mathematics}, pp. 1--42.
  American Mathematical Society, 2011.
  doi:10.1090/conm/558/11047.

\bibitem{ar:BaloghBollobasWeinreich2000}
Balogh J, Bollob{\'a}s B, Weinreich D.
\newblock The Speed of Hereditary Properties of Graphs.
\newblock \emph{J. Comb. Theory, Ser. B}, 2000.
\newblock \textbf{79}(2):131--156.
doi:10.1006/jctb.2000.1952.

\bibitem{ar:BaloghBollobasWeinreich2001}
Balogh J, Bollob{\'a}s B, Weinreich D.
\newblock The Penultimate Rate of Growth for Graph Properties.
\newblock \emph{Eur. J. Comb.}, 2001.
\newblock \textbf{22}(3):277--289.
doi:10.1006/eujc.2000.0476.

\bibitem{balogh2002measures}
Balogh J, Bollob{\'a}s B, Weinreich D.
\newblock Measures on monotone properties of graphs.
\newblock \emph{Discrete Applied Mathematics}, 2002.
\newblock \textbf{116}(1-2):17--36.
doi:10.1016/S0166-218X(01)00175-5.

\bibitem{scheinerman1994size}
Scheinerman ER, Zito J.
\newblock On the size of hereditary classes of graphs.
\newblock \emph{Journal of Combinatorial Theory, Series B}, 1994.
\newblock \textbf{61}(1):16--39.
doi:10.1006/jctb.1994.1027.

\bibitem{ar:Redfield27}
Redfield J.
\newblock The theory of group-reduced distributions.
\newblock \emph{Amer. J. Math.}, 1927.
\newblock \textbf{49}:433--455.

\bibitem{ar:Read59}
Read R.
\newblock The enumeration of locally restricted graphs, I.
\newblock \emph{J. London Math. Soc.}, 1959.
\newblock \textbf{34}:417--436.

\bibitem{ar:Read60}
Read R.
\newblock The enumeration of locally restricted graphs, II.
\newblock \emph{J. London Math. Soc.}, 1960.
\newblock \textbf{35}:344--351.

\bibitem{bk:PolyaRead}
Polya G, Read R.
\newblock Combinatorial enumeration of groups, graphs, and chemical compounds.
\newblock Springer, 1987.
doi:10.1007/978-1-4612-4664-0.

\bibitem{bk:HararyPalmer}
Harary F, Palmer E.
\newblock Graphical Enumeration.
\newblock Academic Press, 1973.
ISBN:978-0-12-324245-7.

\bibitem{ar:Fischer02}
Fischer E.
\newblock The {S}pecker-{B}latter theorem does not hold for quaternary  relations.
\newblock \emph{Journal of Combinatorial Theory, Series A}, 2003.
\newblock \textbf{103}(1):121--136.
doi:10.1016/S0097-3165(03)00075-X.

\bibitem{ar:Specker2005}
Specker E.
\newblock Modular Counting and Substitution of Structures.
\newblock \emph{Combinatorics, Probability and Computing}, 2005.
\newblock \textbf{14}:203--210.
doi:10.1017/S0963548304006698.

\bibitem{pr:FischerMakowsky03}
Fischer E, Makowsky JA.
\newblock The Specker-Blatter Theorem revisited.
\newblock In: COCOON, volume 2697 of \emph{Lecture Notes in Computer Science}.
  Springer, 2003 pp. 90--101.
  doi:10.1007/3-540-45071-8\_11.

\bibitem{bk:Libkin2004}
Libkin L.
\newblock Elements of Finite Model Theory.
\newblock Springer, 2004.
doi:10.1007/978-3-662-07003-1.

\bibitem{alikhani2009introduction}
Alikhani S, Peng Yh.
\newblock Introduction to domination polynomial of a graph.
\newblock \emph{arXiv preprint arXiv:0905.2251}, 2009.

\bibitem{oboudi2016roots}
Oboudi MR.
\newblock On the roots of domination polynomial of graphs.
\newblock \emph{Discrete Applied Mathematics}, 2016.
\newblock \textbf{205}:126--131.
doi:10.1016/j.dam.2015.12.010.

\bibitem{arratia2004interlace}
Arratia R, Bollob{\'a}s B, Sorkin GB.
\newblock The interlace polynomial of a graph.
\newblock \emph{Journal of Combinatorial Theory, Series B}, 2004.
\newblock \textbf{92}(2):199--233.
doi:10.1016/j.jctb.2004.03.003.

\bibitem{courcelle2008multivariate}
Courcelle B.
\newblock A multivariate interlace polynomial and its computation for graphs of
  bounded clique-width.
\newblock \emph{The Electronic Journal of Combinatorics}, 2008. \textbf{15}
\newblock pp. R69--R69.
 doi:10.37236/793.

\bibitem{fischer2003specker}
Fischer E, Makowsky JA.
\newblock The Specker-Blatter Theorem Revisited.
\newblock In: Computing and Combinatorics: 9th Annual International Conference,
  COCOON 2003, volume 2697 of \emph{Springer Lecture Notes in Computer
  Science}. 2003 pp. 90--101.    doi:10.1007/3-540-45071-8\_11.
\end{thebibliography}
\end{document}